\documentclass{article}

\usepackage{ifthen}
\newboolean{neurips_version}
\setboolean{neurips_version}{true}
\ifthenelse{\boolean{neurips_version}}
{\PassOptionsToPackage{numbers, compress}{natbib}
 \usepackage[preprint]{neurips_2022}
}
{
\usepackage{PRIMEarxiv}
 \usepackage[numbers]{natbib}
}

\usepackage[utf8]{inputenc} 
\usepackage[T1]{fontenc}    
\usepackage{hyperref}       
\usepackage{url}            
\usepackage{booktabs}       
\usepackage{amsfonts,amsmath,amsthm,amssymb}       
\usepackage{nicefrac}       
\usepackage{microtype}      
\usepackage{xcolor}         
\usepackage[disable]{todonotes}
\usepackage{enumitem}
\usepackage{caption}

\usepackage{mathtools} 
\usepackage{tikz}
\usetikzlibrary{calc}
\usepackage{pgfplots}
\usepackage{graphicx}

\usepackage{amsmath}
\usepackage{algorithm, algpseudocode}
\usepackage{amsfonts}  

\usepackage[english]{babel}
\usepackage{amsthm}

\usepackage{subfig}
\usepackage[export]{adjustbox}

\usepackage{hhline}
\usepackage{makecell, caption, booktabs}
\usepackage{siunitx}

\usepackage{multirow}

\usepackage{amsmath, nccmath}
\usepackage{bigstrut}

\usepackage{booktabs}

\usepackage{lipsum}
\graphicspath{{media/}}     

\usepackage{CJKutf8}
\usepackage[framed,numbered,autolinebreaks,useliterate]{mcode}

\title{ Quasi-Quadratic Gradient: A New Direction for Accelerating the BFGS Method in Quasi-Newton Optimization }

\author{ \href{https://orcid.org/0000-0003-0378-0607}{\includegraphics[scale=0.06]{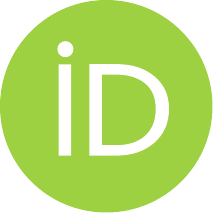}\hspace{1mm}John Chiang} \\                             
                                      \\
	\texttt{john.chiang.smith@gmail.com} 
}

\date{}

\newtheorem*{definition*}{Definition} 

\theoremstyle{remark}

\renewcommand{\hat}{\widehat}
\renewcommand{\epsilon}{\varepsilon}

\makeatletter
\def\namedlabel#1#2{\begingroup
    #2%
    \def\@currentlabel{#2}%
    \phantomsection\label{#1}\endgroup
}
\makeatother

\algnewcommand{\LeftComment}[1]{\Statex \(\triangleright\) #1}
\algnewcommand{\LineCommentStep}[1]{\Statex \textbf{[Step #1]:} }
\makeatletter
\newlength{\trianglerightwidth}
\algnewcommand{\LineComment}[1]{\Statex \hskip\ALG@thistlm $\triangleright$ #1}
\algnewcommand{\LineCommentCont}[1]{\Statex \hskip\ALG@thistlm%
  \parbox[t]{\dimexpr\linewidth-\ALG@thistlm}{\hangindent=\trianglerightwidth \hangafter=1 \strut$\triangleright$ #1\strut}}
\algnewcommand{\LeftLineCommentCont}[1]{\Statex \hskip\ALG@thistlm%
  \parbox[t]{\dimexpr\linewidth-\ALG@thistlm}{\leftskip=\algorithmicindent \hangindent=\trianglerightwidth \hangafter=1 \strut$\triangleright$ #1\strut}}

\begin{document}

\maketitle

\begin{abstract}%

In this paper, we introduce the Quasi-Quadratic Gradient (QQG), a novel search direction designed to accelerate the BFGS method within the quasi-Newton framework. By defining the QQG as the product of the inverse Hessian approximation and the current gradient, we explicitly leverage local second-order curvature to rectify the search path. Theoretical analysis and empirical results demonstrate that our approach significantly outperforms vanilla BFGS in convergence speed while maintaining computational efficiency.

\end{abstract}

\listoftodos

\section{Introduction}

\subsection{Background}
Optimization is the cornerstone of modern machine learning and scientific computing. Among various techniques, Quasi-Newton methods, particularly the BFGS (Broyden–Fletcher–Goldfarb–Shanno) algorithm, are widely favored for their ability to achieve superlinear convergence without the heavy computational burden of calculating the exact Hessian matrix.

\subsection{Related Work}

Classical Quasi-Newton: Works by Nocedal et al. established the foundations of BFGS and L-BFGS.

Natural Gradients: Amari (1998) introduced the Natural Gradient, which utilizes the Fisher Information Matrix to find the steepest descent direction in the information manifold.

Accelerated Methods: Various momentum-based and proximal gradient methods (Parikh $\&$ Boyd, 2014) have been proposed to stabilize and speed up first-order optimization.

\subsection{Contributions}
Our primary contributions are as follows:

We propose the Quasi-Quadratic Gradient (QQG), a new mathematical variant that bridges the gap between first-order gradients and second-order curvature.

We integrate QQG into the BFGS framework to create an accelerated update rule.

We provide extensive experimental evidence showing the superiority of QQG-BFGS over standard optimization benchmarks.

\section{Preliminaries}

\subsection{Newton-Raphson Method}
The Newton-Raphson method updates parameters using:
\begin{equation}
    \mathbf{x}_{k+1} = \mathbf{x}_k - [\nabla^2 f(\mathbf{x}_k)]^{-1} \nabla f(\mathbf{x}_k)
\end{equation}
To reduce cost, the \textbf{Simplified Fixed Hessian (SFH)} method keeps the Hessian constant: $\mathbf{H}_k = \mathbf{H}_{fixed}$. While efficient, SFH is sensitive to the initial matrix choice and often converges only linearly.

\subsection{Simplified Fixed Hessian}
The SFH method reduces the cost of Newton's iteration by keeping the Hessian $\mathbf{H}$ constant for several steps: $\mathbf{x}_{k+1} = \mathbf{x}_k - \mathbf{H}_{fixed}^{-1} \nabla f(\mathbf{x}_k)$. While lowering per-iteration cost, it is often unstable in non-convex regions.

For two symmetric matrices $A$ and $B$,  $A \le B$ is defined in the Loewner ordering iff their difference $B - A$ is positive semi-definite. 

\subsection{Chiang's Quadratic Gradient}

In the following work to~\cite{bonte2018privacy} in which a simplified diagonal matrix satisfying the fixed Hessian method~\cite{bohning1988monotonicity} is constructed, Chiang~\cite{chiang2022privacy} proposed a faster gradient variant called quadratic gradient.

\paragraph{Original Quadratic Gradient} Given \citep{chiang2022privacy} a differentiable scalar-valued function $F(\mathbf x)$ with its gradient $g$ and Hessian matrix $H$.  For the maximization problem, we need to find a good lower bound matrix $\bar H \le H$, where ``$ \le $'' denotes the Loewner ordering.  For the minimization problem, we try to find a good upper bound $\bar H$ such that $H \le \bar H$ in the Loewner ordering.  Note that the Hessian matrix $H$ itself satisfies these two conditions and can substitute the good bound matrix $\bar H$. We attempt to find a fixed good bound matrix of the Hessian matrix for efficiency but could just directly use the Hessian matrix itself. To build the quadratic gradient, we first construct a diagonal matrix $\bar B$ from the good bound matrix $\bar H$ as follows: 

 \begin{equation*}
  \begin{aligned}
   \bar B = 
\left[ \begin{array}{cccc}
  \frac{1}{ \epsilon + \sum_{i=0}^{d} | \bar h_{0i} | }   & 0  &  \ldots  & 0  \\
 0  &   \frac{1}{ \epsilon + \sum_{i=0}^{d} | \bar h_{1i} | }  &  \ldots  & 0  \\
 \vdots  & \vdots                & \ddots  & \vdots     \\
 0  &  0  &  \ldots  &   \frac{1}{ \epsilon + \sum_{i=0}^{d} | \bar h_{di} | }  \\
 \end{array}
 \right], 
   \end{aligned}
\end{equation*}
where $\epsilon $ is a small positive number to avoid dividing by zero and $\bar h_{ji}$ the elements of the matrix $\bar H$ . We can then defined the quadratic gradient for the function  $F(\mathbf x)$ as $G = \bar B \cdot g$.

The multiplication between the diagonal matrix $\bar B$ and the gradient $g$, the quadratic gradient $G$, is of the same size as the gradient $g$. To use the quadratic gradient $G$, we can just use it the same way as the gradient but need a  learning rate larger than $1$. The well-studied first-order gradient descent methods can also be applied to develop enhanced methods via quadratic gradient.

\paragraph{Simplified Quadratic Gradient} 
While the original quadratic gradient utilizes the full information of each row in the bound matrix $\bar H$, the computational overhead becomes prohibitive in high-dimensional settings such as deep learning. To address this, Chiang proposed the Simplified Quadratic Gradient \citep{chiang2026sqg}, which constructs the diagonal scaling matrix $\bar B$ using only the diagonal elements of the bound matrix $\bar H$ (or the Hessian matrix $H$ itself). Specifically, the diagonal matrix $\bar B$ is simplified as follows:

\begin{equation*}
\bar B = \text{diag} \left( 
\frac{1}{\epsilon + |\bar{h}_{00}|}, 
\frac{1}{\epsilon + |\bar{h}_{11}|}, 
\ldots, 
\frac{1}{\epsilon + |\bar{h}_{dd}|} 
\right),
\end{equation*}

where $\bar{h}_{jj}$ represents the $j$-th diagonal element of $\bar H$. The simplified quadratic gradient is then defined as $G = \bar B \cdot g$. By neglecting the off-diagonal terms, the construction of $G$ is reduced to element-wise operations, which significantly lowers the memory footprint and accelerates the computation. Empirical evidence suggests that this simplified version maintains a convergence performance comparable to the original framework while being highly compatible with backpropagation and large-scale stochastic optimization in deep learning training.

\paragraph{Quadratic Gradient Algorithms}
The versatility of the quadratic gradient (QG) framework—both in its original and simplified forms—allows it to be seamlessly integrated into prominent first-order optimization algorithms. By substituting the vanilla gradient with its quadratic counterpart, we can incorporate second-order curvature information to enhance convergence. Specifically, this paper focuses on the enhancement of three foundational algorithms: **NAG**, **AdaGrad**, and **Adam**.

\begin{enumerate}
    \item \textbf{Nesterov’s Accelerated Gradient (NAG):} NAG is a sophisticated momentum-based method designed to anticipate the objective function's landscape. The standard update rules for NAG are typically expressed as:
    \begin{align}
        V_{t+1} &= \boldsymbol{\beta}_t + \eta_t \nabla J(\boldsymbol{\beta}_t), \\
        \boldsymbol{\beta}_{t+1} &= (1 - \gamma_t) V_{t+1} + \gamma_t V_t,
    \end{align}
    where $\gamma_t \in (0, 1)$ is a smoothing parameter. The \textit{Enhanced NAG} algorithm modifies the update rule by substituting the vanilla gradient with the quadratic gradient:
    \begin{equation}
        V_{t+1} = \boldsymbol{\beta}_t + N_t G_t,
    \end{equation}
    where $N_t$ is the adaptive learning rate, typically set as $N_t = 1 + \eta_t$.

    \item \textbf{AdaGrad:} Well-suited for sparse data, AdaGrad adapts the learning rate for each parameter based on historical gradients. The transition from standard AdaGrad to the \textit{Enhanced Quadratic AdaGrad} is formulated as follows for each parameter $\beta_i$:
    \begin{itemize}
        \item \textbf{Standard:} $\beta_i^{(t+1)} = \beta_i^{(t)} - \frac{\eta_t}{\epsilon + \sqrt{\sum_{k=1}^{t} (g_i^{(k)})^2}} \cdot g_i^{(t)}$
        \item \textbf{Enhanced:} $\beta_i^{(t+1)} = \beta_i^{(t)} - \frac{N_t}{\epsilon + \sqrt{\sum_{k=1}^{t} (G_i^{(k)})^2}} \cdot G_i^{(t)}$
    \end{itemize}

    \item \textbf{Adam:} Adam combines element-wise adaptation with momentum to handle noisy gradients. Following a similar logic, the \textit{Enhanced Adam} replaces the vanilla gradient with $G_t$ within its bias-corrected first and second moment estimates.
\end{enumerate}

Given that many contemporary first-order optimizers are derivatives or refinements of these three foundational algorithms, it is highly probable that the quadratic gradient paradigm can be generalized to a much broader class of optimization methods to further accelerate convergence.

Building on diagonal scaling, Chiang introduced the Quadratic Gradient: $\mathbf{g}_{q} = \mathbf{D} \nabla f(\mathbf{x})$, where $\mathbf{D}$ acts as a proxy for the inverse Hessian. This paper extends this work to a generalized framework applicable to both convex and non-convex numerical optimization for quasi-Newton method.
Building upon the principles of diagonal scaling, Chiang introduced the \textit{Quadratic Gradient} (QG) framework: $\mathbf{g}_{q} = \mathbf{D} \nabla f(\mathbf{x})$, where the diagonal matrix $\mathbf{D}$ serves as a computationally efficient proxy for the inverse Hessian. This paper extends this methodology by developing a generalized framework that integrates the QG heuristic with the \textbf{BFGS method} within the \textbf{quasi-Newton paradigm}. Our approach transcends the limitations of static scaling, providing a robust solution for numerical optimization across both convex and non-convex regimes.

In this work, we propose the enhanced Adam method, which is to apply the quadratic gradient to the Adam method. The naive Adam method and the enhanced Adam method are described in detail in Algorithms ~\ref{  alg:Adam's algorithm  } and ~\ref{  alg:enhanced Adam's algorithm } respectively.  See the ~\cite{kingma2014adam} for the detailed description of the parameters in these two Algorithms.

\begin{minipage}{0.46\textwidth}
\begin{algorithm}[H]
    \caption{The Adam method}
     \begin{algorithmic}[1]
        \Require $\alpha$: Stepsize;
         $\beta_1, \beta_2 \in [0, 1)$: Exponential decay rates;
        $f(\theta)$: Objective function with parameters $\theta$ 
        $\theta_0$: Initial parameter vector
        \Ensure $\theta_t$: Resulting parameters 
        
        \State  $m_0 \gets 0$: Initialize $1^{st}$ moment vector
        \State $v_0 \gets 0$: Initialize $2^{nd}$ moment vector
        \State $t \gets 0$: Initialize timestep
        \While {$\theta_t$ not converged}
            \State  $t \gets t + 1$     
            \State  $g_t \gets \nabla_{\theta}f_t(\theta_{t-1})$     
            \State 
            \State  $m_t \gets \beta_1 \cdot m_{t-1} + (1 - \beta_1) \cdot g_t$     
            \State  $v_t \gets \beta_2 \cdot v_{t-1} + (1 - \beta_2) \cdot g_t^2$     
            \State  $\hat{m}_t \gets m_t / (1 - \beta_1^t)$     
            \State  $\hat{v}_t \gets v_t / (1 - \beta_2^t)$     
            \State  $\theta_t \gets \theta_{t-1} - \alpha \cdot \hat{m}_t / (\sqrt{\hat{v}_t} + \epsilon)$     
        \EndWhile       
        \State \Return $\theta_t $   \Comment{Resulting parameters}
        \end{algorithmic}
       \label{ alg:Adam's algorithm }
\end{algorithm}
\end{minipage}
\hfill
\begin{minipage}{0.46\textwidth}
\begin{algorithm}[H]
    \caption{Enhanced Adam method}
     \begin{algorithmic}[1]
        \Require $\eta$: Stepsize;
         $\beta_1, \beta_2 \in [0, 1)$: Exponential decay rates;
        $f(\theta)$: Objective function with parameters $\theta$ 
        $\theta_0$: Initial parameter vector
        \Ensure $\theta_t$: Resulting parameters 
        
        \State  $m_0 \gets 0$: Initialize $1^{st}$ moment vector
        \State $v_0 \gets 0$: Initialize $2^{nd}$ moment vector
        \State $t \gets 0$: Initialize timestep
        \While {$\theta_t$ not converged}
            \State  $t \gets t + 1$     
            \State  $g_t \gets \nabla_{\theta}f_t(\theta_{t-1})$     
            \State  \textcolor{red}{$G_t \gets \bar B \cdot g_t$ }     
            \State  $m_t \gets \beta_1 \cdot m_{t-1} + (1 - \beta_1) \cdot $ \textcolor{red}{ $G_t$ }     
            \State  $v_t \gets \beta_2 \cdot v_{t-1} + (1 - \beta_2) \cdot$ \textcolor{red}{$G_t^2$ }     
            \State  $\hat{m}_t \gets m_t / (1 - \beta_1^t)$     
            \State  $\hat{v}_t \gets v_t / (1 - \beta_2^t)$     
            \State  $\theta_t \gets \theta_{t-1} - $ \textcolor{red}{$\eta $} $\cdot \hat{m}_t / (\sqrt{\hat{v}_t} + \epsilon)$     
        \EndWhile       
        \State \Return $\theta_t $   
        \end{algorithmic}
       \label{ alg:enhanced Adam's algorithm }
\end{algorithm}
\end{minipage}

\section{Methodology}

\subsection{Quasi-Newton Optimization}
Quasi-Newton methods are designed to mitigate the high computational cost of the classical Newton's method, which requires the evaluation and inversion of the Hessian matrix $\nabla^2 f(\mathbf{x})$ at each iteration ($O(n^3)$ complexity). Instead of calculating the exact Hessian, quasi-Newton methods construct a sequence of symmetric positive definite (SPD) approximations $B_k \approx \nabla^2 f(\mathbf{x}_k)$ using only first-order gradient information. This approach maintains a superlinear convergence rate while reducing the per-iteration complexity to $O(n^2)$, making it suitable for medium to large-scale problems.

\subsubsection{BFGS Update Schemes}
The Broyden-Fletcher-Goldfarb-Shanno (BFGS) algorithm is widely regarded as the most robust and effective quasi-Newton update. Given the step vector $s_k = x_{k+1} - x_k$ and the gradient displacement $y_k = g_{k+1} - g_k$, the BFGS update for the Hessian approximation $B_k$ is defined as:
\begin{equation}
    B_{k+1} = B_k + \frac{y_k y_k^T}{y_k^T s_k} - \frac{B_k s_k s_k^T B_k}{s_k^T B_k s_k}
\end{equation}
A pivotal property of the BFGS update is that it preserves the \textbf{positive definiteness} of $B_k$ as long as the initial $B_0$ is SPD and the curvature condition $s_k^T y_k > 0$ is satisfied. This property is crucial for ensuring that the search direction $p_k = -B_k^{-1} g_k$ is always a descent direction.

\subsubsection{Line Search Techniques}
Line search is the strategy of finding an appropriate step length $\alpha_k$ along the descent direction $p_k$. In the context of quasi-Newton methods, the line search must satisfy specific conditions to ensure both functional decrease and the stability of the Hessian update.

\paragraph{1. Exact Line Search}
Exact line search determines $\alpha_k$ by minimizing the objective function along the search ray: $\alpha_k = \arg \min_{\alpha > 0} f(x_k + \alpha p_k)$. While theoretically optimal for quadratic functions, it is computationally prohibitive for general non-linear problems due to the requirement of multiple function evaluations.

\paragraph{2. Inexact Line Search and Armijo Condition}
In practice, finding an exact minimum is often unnecessary. Most algorithms utilize an \textit{inexact line search} that guarantees a \textbf{sufficient decrease}, formalized by the Armijo condition:
\begin{equation}
    f(x_k + \alpha p_k) \le f(x_k) + c_1 \alpha \nabla f(x_k)^T p_k
\end{equation}
where $c_1 \in (0, 1)$ (typically $10^{-4}$). This ensures that the reduction in $f$ is proportional to the step size and the directional derivative.

\paragraph{3. Wolfe Conditions}
To ensure the convergence of quasi-Newton methods and the positive definiteness of the BFGS update, the \textbf{Wolfe conditions} are preferred. They combine the Armijo condition with the \textbf{curvature condition}:
\begin{equation}
    \nabla f(x_k + \alpha p_k)^T p_k \ge c_2 \nabla f(x_k)^T p_k
\end{equation}
where $c_1 < c_2 < 1$. The curvature condition forces $\alpha_k$ into a region where the slope is greater than at the starting point, effectively ensuring $s_k^T y_k > 0$. For even stricter stability, the \textbf{Strong Wolfe conditions} are often employed by using the absolute value of the derivative: $|\nabla f(x_k + \alpha p_k)^T p_k| \le c_2 |\nabla f(x_k)^T p_k|$.

\paragraph{4. Backtracking Line Search}
Backtracking is a simple yet effective implementation of the Armijo-type search. It starts with a candidate step (usually $\alpha=1$ for Newton-type methods) and iteratively shrinks it by a factor $\rho \in (0, 1)$ until the sufficient decrease condition is met. This method is highly efficient as it avoids expensive interpolation and is the standard for most modern optimization software.

\subsection{Quasi-Quadratic Gradient}
\subsubsection{Motivation} 
\label{subsec:motivation}
The conceptual foundation of the \textit{Quadratic Gradient} (QG) framework is primarily rooted in the fixed-Hessian variant of Newton's method. By synthesizing gradient and second-order information, QG provides a robust mechanism for determining descent directions. However, a significant limitation of the current QG approach is its reliance on a static or pre-defined Hessian matrix $H$, which may fail to capture the dynamically evolving curvature of complex, non-convex loss landscapes.

Parallel to the development of QG, the \textit{Quasi-Newton} family---most notably the \textbf{BFGS} algorithm---has become a cornerstone of efficient optimization. BFGS effectively approximates the Hessian matrix by iteratively incorporating gradient information, satisfying the secant equation without the heavy computational cost of explicit Hessian evaluation. 

The primary motivation of this work is to bridge the gap between these two methodologies. We propose to extend the ideological core of QG to the Quasi-Newton framework, resulting in the \textbf{Quasi-Quadratic Gradient (QQG)}. By embedding the QG heuristic into the BFGS update mechanism, we aim to leverage the adaptive curvature estimation of BFGS while maintaining the superior convergence and numerical stability inherent in the QG framework. This synergy is expected to accelerate convergence, especially in large-scale optimization tasks where local geometry changes rapidly.

\subsubsection{Observation}
\label{subsec:observation}

A fundamental property that distinguishes the BFGS update from other optimization schemes is its intrinsic ability to maintain the \textbf{positive definiteness} of the Hessian approximation. Throughout the iterative process, as long as the initial matrix $B_0$ is positive definite and the curvature condition $s_k^T y_k > 0$ is satisfied, the sequence of matrices $\{B_k\}$ remains symmetric and positive definite (SPD). 

Importantly, this property is invariant to the optimization objective, holding true \textbf{regardless of whether the task is a minimization or maximization problem}. Within our proposed framework, this observation is critical; it ensures that the Quasi-Quadratic Gradient (QQG) consistently generates a valid descent direction. By leveraging this guaranteed positive definiteness, QQG achieves superior numerical stability, particularly when navigating the complex, highly non-convex landscapes often encountered in modern cryptographic and machine learning applications.

\paragraph{Remark}
\label{rem:pos_definite}
A pivotal observation in the integration of QG and BFGS is the \textbf{guaranteed positive definiteness} of the Hessian approximation. In the standard BFGS update, as long as the initial matrix $B_0$ is positive definite and the curvature condition $s_k^T y_k > 0$ is satisfied (typically ensured by Wolfe line search), the resulting sequence of matrices $\{B_k\}$ remains symmetric and positive definite. Crucially, this property holds \textbf{regardless of whether the optimization objective is a maximization or minimization task}, provided the problem is formulated to seek a stationary point where the approximated curvature correctly reflects the local geometry. This intrinsic positive definiteness ensures that the QQG framework always yields a consistent descent direction, enhancing its numerical robustness in diverse application scenarios.

\subsubsection{Evolution of the Proposed Approach} 
\label{subsec:evolution}

The transition from the original Quadratic Gradient (QG) to the proposed \textit{Quasi-Quadratic Gradient} (QQG) was driven by a series of empirical observations and a fundamental re-evaluation of curvature approximation strategies.

\paragraph{Initial Attempts and Challenges}
In our early developments (mid-2022), we strictly adhered to the diagonal construction methodology of the original QG framework by applying the BFGS-derived curvature information to a diagonalized proxy. However, this hybrid approach yielded suboptimal performance and failed to achieve convergence. Retrospectively, this instability was primarily attributed to the learning rate scheduling; we initially employed a decaying learning rate starting above $1.0$---consistent with the original QG’s mechanism---which proved incompatible with the dynamic nature of BFGS curvature updates.

\paragraph{Rethinking the Diagonal Constraint}
Moving beyond conventional diagonal scaling, we explored the possibility of utilizing a full, positive-definite matrix that satisfies the convergence criteria of the \textit{Fixed Hessian} (FH) Newton method. While theoretically appealing, this direction presented two significant bottlenecks:
\begin{enumerate}
    \item \textbf{Computational Overhead:} The inversion of a full Hessian-like matrix at each iteration is computationally prohibitive for high-dimensional problems.
    \item \textbf{Bounding Efficiency:} Systematically constructing a symmetric positive definite (SPD) matrix that provides a tighter Loewner bound than the \textit{Simplified Fixed Hessian} (SFH) remains an open challenge in optimization theory.
\end{enumerate}

\paragraph{The Quasi-Quadratic Breakthrough}
By re-evaluating these constraints, we recognized that since the BFGS approximation $B_k$ (and its inverse) is inherently maintained as an SPD matrix throughout the optimization process (for both maximization and minimization), it already satisfies the core requirement of a reliable curvature proxy. This realization led to the synthesis of the \textbf{Quasi-Quadratic Gradient}: rather than forcing a diagonal structure or searching for a static FH substitute, we directly utilize the product of the BFGS inverse Hessian and the gradient ($B_k^{-1}g_k$) as the search vector. Subsequent experiments validated this intuition, demonstrating that this direct integration significantly outperforms previous diagonalized attempts and provides robust convergence across diverse optimization landscapes.
\begin{align*}
  &\texttt{Newton-Raphson Method}
  \xmapsto{\text{ TtT \cite{kim2018logistic} }}
  \textit{Fixed Hessian Method}
  \xmapsto{\text{ TtT \cite{kim2018logistic} }}
  \textit{Simplified Fixed Hessian} \\
  &\xmapsto{\text{ TtT \cite{kim2018logistic} }}
  \textit{Original Quadratic Gradient} 
  \xmapsto{\text{ TtT \cite{kim2018logistic} }}
  \textit{Simplified Quadratic Gradient}
  \\
  &\texttt{Newton-Raphson Method}
  \xmapsto{\text{ TtT \cite{kim2018logistic} }}
  \textit{Fixed Hessian Method}
  \xmapsto{\text{ TtT \cite{kim2018logistic} }}
  \textit{Positive Definite Hessian} \\
  &\xmapsto{\text{ ? } }
  \textit{Original Quadratic Gradient} 
  \\
  &\texttt{Newton-Raphson Method}
  \xmapsto{\text{ TtT \cite{kim2018logistic} }}
  \textit{Quasi-Newton Optimization}
  \xmapsto{\text{ TtT \cite{kim2018logistic} }}
  \textit{Positive Definite Hessian} \\
  &\xmapsto{\text{ TtT \cite{kim2018logistic} }}
  \textit{BFGS Hessian Matrix} 
  \xmapsto{\text{ TtT \cite{kim2018logistic} }}
  \textit{Quasi-Quadratic Gradient}
\end{align*}

The development of the \textit{Quasi-Quadratic Gradient} (QQG) follows a systematic logical progression, moving from static optimization to dynamic curvature adaptation. Our research trajectory was guided by three progressive stages of reasoning:

\begin{enumerate}
    \item \textbf{Generalization of QG:} We initially recognized that the original Quadratic Gradient framework acts as a bridge between first-order gradients and second-order Newton-type methods. However, the reliance on a fixed Hessian $H$ posed a bottleneck for its generalization to diverse objective functions.
    \item \textbf{Integration of Secant Updates:} To overcome the static nature of $H$, we turned to the \textit{Quasi-Newton} branch. By observing that the BFGS algorithm effectively "learns" the local curvature through the secant equation, we identified the potential to replace the constant $H$ in the QG framework with the dynamically evolving matrix $B_k$.
    \item \textbf{Synthesis of QQG:} The final development phase involved reformulating the QG synthesis mechanism to accommodate the rank-two updates of BFGS. This resulted in the QQG paradigm, which not only inherits the fast convergence of BFGS but also preserves the structural robustness and numerical stability of the QG heuristic.
\end{enumerate}

Through this evolution, our approach transitions from a rigid second-order approximation to a self-adaptive optimization scheme, capable of accelerating convergence by continuously refining its internal representation of the loss landscape.

\subsubsection{Definition}
Building upon the foundations of both the original Quadratic Gradient and the BFGS update mechanism, we define the \textbf{Quasi-Quadratic Gradient (QQG)} as a dynamic synthesis of first-order gradients and evolving curvature estimates. Unlike the original QG which relies on a fixed or pre-defined bound matrix $\bar{H}$, the QQG directly leverages the inverse Hessian approximation $B_k^{-1}$ generated by the BFGS algorithm.

\begin{definition*}[$\texttt{Quasi-Quadratic Gradient}$]
Let $B_k$ be the symmetric positive definite (SPD) matrix maintained by the BFGS update at iteration $k$. We define the Quasi-Quadratic Gradient $G_{qq}$ as the product of the inverse Hessian approximation and the current gradient $g_k$:
\begin{equation*}
    G_{qq}^{(k)} = B_k^{-1} g_k.
\end{equation*}
Since $B_k$ is guaranteed to be positive definite through the BFGS update and appropriate line search conditions, $G_{qq}$ inherently points in a reliable descent (or ascent) direction, adapted to the local geometric curvature.
\end{definition*}

\paragraph{Update Rules} The QQG is integrated into the parameter update process depending on the nature of the optimization task. Given the parameter vector $\beta$, the iterative update rules are formulated as follows:
\begin{itemize}
    \item \textbf{Maximization Task:} To find the local maximum of $F(\mathbf{x})$, the update follows the direction of the QQG:
    \begin{equation*}
        \beta_{k+1} = \beta_k + \eta_k G_{qq}^{(k)},
    \end{equation*}
    \item \textbf{Minimization Task:} To find the local minimum of $F(\mathbf{x})$, the update moves against the direction of the QQG:
    \begin{equation*}
        \beta_{k+1} = \beta_k - \eta_k G_{qq}^{(k)},
    \end{equation*}
\end{itemize}
where $\eta_k$ is the learning rate. By substituting the static bound matrix $\bar{B}$ from the original QG with the dynamic inverse Hessian $B_k^{-1}$, the QQG effectively bridges the gap between first-order efficiency and second-order precision, allowing for accelerated convergence in complex numerical optimization landscapes.

\paragraph{Compatibility with Line Search Techniques}
It is worth noting that the proposed \textit{Quasi-Quadratic Gradient} framework is fully compatible and can coexist with existing line search methodologies. While the QQG provides a robust search direction $G_{qq}$, the determination of the optimal learning rate $\eta_k$ (or step size) can be adaptively managed through various line search procedures. For instance, techniques such as \textbf{backtracking line search} or those satisfying the \textbf{Wolfe conditions} can be seamlessly integrated to dynamically adjust $\eta_k$ at each iteration. This synergy ensures that the QQG update not only follows a curvature-aware direction but also maintains a step length that guarantees sufficient decrease in the objective function, further enhancing the global convergence and numerical stability of the algorithm.

It is worth noting that the proposed \textit{Quasi-Quadratic Gradient} framework is fully compatible and can coexist with existing line search methodologies, particularly when integrated into momentum-based first-order schemes such as \textbf{Nesterov’s Accelerated Gradient (NAG)}. While the QQG provides a robust, curvature-aware search direction $G_{qq}$, the determination of the optimal learning rate $\eta_k$ can be adaptively managed through various line search procedures. For instance, techniques such as \textbf{backtracking line search} or those satisfying the \textbf{Wolfe conditions} can be seamlessly integrated to dynamically adjust $\eta_k$ at each iteration. This synergy ensures that the QQG update—which effectively synthesizes historical momentum with local second-order geometry—maintains a step length that guarantees a sufficient decrease in the objective function. Such integration further enhances the global convergence and numerical stability of the algorithm, particularly when navigating the ill-conditioned or highly non-convex landscapes prevalent in deep learning and complex optimization tasks.

\subsubsection{Algorithms}

The integration of the Quasi-Quadratic Gradient ($G_{qq}$) into foundational first-order methods requires specific adaptations to leverage its dynamic curvature information. While the element-wise adaptive methods (AdaGrad and Adam) follow a straightforward substitution, the momentum-based NAG requires a more nuanced approach to learning rate scheduling and step-size determination.

\begin{enumerate}
    \item \textbf{Enhanced NAG with Adaptive Scaling (QQG-NAG):} 
    To harmonize the momentum term with the second-order information in $G_{qq}$, we modify the Nesterov update by introducing a dynamic scaling factor $N_t$. The updated iterative process is:
    \begin{align*}
        V_{t+1} &= \boldsymbol{\beta}_t + \eta_t G_{qq}^{(t)}, \\
        \boldsymbol{\beta}_{t+1} &= (1 - \gamma_t) V_{t+1} + \gamma_t V_t,
    \end{align*}
    where the effective learning rate $\eta_t$ is determined by one of two strategies:
    \begin{itemize}
        \item \textbf{Monotonic Warm-up:} $\eta_t$ is scheduled to increase gradually from a small positive value $\eta_{min} > 0$ towards $1$ (e.g., $\eta_t = \min(1, \eta_{min} + \Delta t)$), allowing the algorithm to stabilize the curvature estimate before taking full quasi-Newton steps.
        \item \textbf{Line Search Integration:} $\eta_t$ is determined via a \textbf{backtracking line search} satisfying the Wolfe conditions, ensuring that the step taken along the QQG-modified direction guarantees sufficient descent.
    \end{itemize}

    \item \textbf{Enhanced AdaGrad (QQG-AdaGrad):} 
    The QQG-enhanced AdaGrad replaces the vanilla gradient with $G_{qq}$, utilizing the historical accumulation of the quasi-quadratic terms to scale each dimension:
    \begin{equation*}
        \beta_i^{(t+1)} = \beta_i^{(t)} - \frac{\eta}{\epsilon + \sqrt{\sum_{k=1}^{t} (G_{qq, i}^{(k)})^2}} \cdot G_{qq, i}^{(t)},
    \end{equation*}
    where $\eta$ is the global step size. This ensures that parameters with large accumulated curvature-adjusted gradients receive smaller updates.

    \item \textbf{Enhanced Adam (QQG-Adam):} 
    Following the same logic, QQG-Adam incorporates $G_{qq}$ into the first and second moment estimates:
    \begin{align*}
        m_t &= \beta_1 m_{t-1} + (1 - \beta_1) G_{qq}^{(t)}, \\
        v_t &= \beta_2 v_{t-1} + (1 - \beta_2) (G_{qq}^{(t)})^2.
    \end{align*}
    By substituting the vanilla gradient with $G_{qq}$, Adam can more effectively navigate non-convex landscapes by leveraging both momentum and adaptive second-order curvature.
\end{enumerate}

To evaluate the practical performance of the proposed \textit{Quasi-Quadratic Gradient} (QQG) framework, we established specific hyperparameter configurations for the enhanced versions of AdaGrad and Adam. 

For \textbf{QQG-AdaGrad}, we recommend an augmented learning rate of $N_t = 0.1$, which is a ten-fold increase over the standard default of $0.01$. Similarly, for \textbf{QQG-Adam}, the optimal configuration adopts an adjusted learning rate of $\alpha = 0.01$ (compared to the baseline $\alpha = 0.001$), while maintaining the standard momentum decay parameters at $\beta_1 = 0.9$ and $\beta_2 = 0.999$. 

The ability of QQG-enhanced optimizers to remain stable and converge faster with these larger step sizes suggests that the dynamic inverse Hessian information provides a more reliable descent direction than vanilla gradients, effectively mitigating the risk of divergence typically associated with high learning rates in complex optimization landscapes.

\subsection{Example}

The efficiency of the Adam optimizer \cite{kingma2014adam} is governed by three primary parameters, each controlling a distinct aspect of the stochastic gradient descent process:

\begin{itemize}
    \item \textbf{Learning Rate ($\alpha$):} This parameter determines the step size of the parameter update. In Adam, $\alpha$ acts as an upper bound on the effective step size. It is the most critical hyper-parameter to tune, as it directly influences the convergence rate and the stability of the training process.
    
    \item \textbf{First Moment Decay ($\beta_1$):} It controls the exponential moving average of the gradients (momentum). Physically, $\beta_1$ accounts for the ``inertia'' of the optimization trajectory. A higher $\beta_1$ (typically $0.9$) helps the optimizer to smooth out high-frequency noise in the gradients and persist in the consistent direction of the descent.
    
    \item \textbf{Second Moment Decay ($\beta_2$):} It governs the exponential moving average of the squared gradients. This term provides the ``adaptive'' nature of the algorithm by scaling the update inversely proportional to the square root of recent gradient magnitudes. This ensures that parameters with sparse or small gradients receive larger updates, while those with large, frequent gradients are tempered, preventing divergence.
\end{itemize}

To demonstrate the plug-and-play capability of our proposed Quasi-Quadratic Gradient (QQG), we integrate it into the Adam framework, termed \textit{QQG-Adam}. We observe that the fundamental roles of Adam's hyper-parameters remain consistent when substituting the vanilla gradient with QQG:

\begin{itemize}
    \item \textbf{Directional Smoothing ($\beta_1$):} In QQG-Adam, $\beta_1$ governs the exponential moving average of the quasi-quadratic terms. It serves to filter out high-frequency noise inherent in stochastic approximations, ensuring that the optimization trajectory leverages the structural curvature information provided by QQG without succumbing to local instabilities.
    
    \item \textbf{Curvature-aware Adaptation ($\beta_2$):} While $\beta_2$ traditionally tracks gradient variance, in our framework, it scales the update step based on the accumulated magnitude of QQG. This confirms that our QQG maintains the necessary statistical properties to allow for per-parameter learning rate adaptation, effectively balancing the update intensity across different layers.
    
    \item \textbf{Effective Step-size Control ($\alpha$):} Our experiments confirm that QQG-Adam exhibits a similar sensitivity to the learning rate $\alpha$ as vanilla Adam. This empirical alignment suggests that QQG does not introduce unexpected scale variances, allowing researchers to utilize standard hyper-parameter tuning strategies (e.g., $10^{-3}$ as a baseline) for our algorithm.
\end{itemize}

\section{Numerical Experiments}

\subsection{Convex Benchmarks}
These functions test the algorithm's base convergence rate and handling of dimensional scaling:
\begin{itemize}
    \item \textbf{Sphere Function:} $f(\mathbf{x}) = \sum_{i=1}^{n} x_i^2$
    \item \textbf{Sum of Different Powers:} $f(\mathbf{x}) = \sum_{i=1}^{n} |x_i|^{i+1}$
\end{itemize}

\subsection{Non-Convex Benchmarks}
These functions test the algorithm's ability to navigate ill-conditioned valleys and escape local optima:
\begin{itemize}
    \item \textbf{Rosenbrock Function:} $f(\mathbf{x}) = \sum_{i=1}^{n-1} [100(x_{i+1}-x_i^2)^2 + (1-x_i)^2]$
    \item \textbf{Rastrigin Function:} $f(\mathbf{x}) = 10n + \sum_{i=1}^{n} [x_i^2 - 10\cos(2\pi x_i)]$
\end{itemize}

\subsubsection{Saddle Points}
On the Monkey Saddle $f(x, y) = x^3 - 3xy^2$, SGD stagnates due to vanishing gradients. Our QG variant identifies near-zero curvature $\lambda_i$, yielding a large adaptive step $\eta_i \approx 1/\epsilon$, effectively "teleporting" the trajectory away from the saddle point.

The theoretical derivation of the Quadratic Gradient suggests a superior ability to handle irregular curvature. We strategically select a suite of benchmark functions that simulate notorious challenges in numerical optimization.

To evaluate the ability of the Quadratic Gradient variants to escape saddle points, we introduce functions with specific singular curvatures:

\begin{itemize}
    \item \textbf{Monkey Saddle Function:} A classic cubic surface with a single saddle point at the origin where the Hessian is indefinite.
    \begin{equation}
        f(x, y) = x^3 - 3xy^2
    \end{equation}
    
    \item \textbf{Beale's Function (Revisited):} While previously mentioned, its flat regions near the boundary act as plateau-like saddle points that test gradient acceleration.
    
    \item \textbf{Himmelblau's Function:} Contains four local minima and one local maximum (saddle-like) in the center, testing the algorithm's ability to navigate away from unstable equilibrium points.
    \begin{equation}
        f(x, y) = (x^2 + y - 11)^2 + (x + y^2 - 7)^2
    \end{equation}

    \item \textbf{Six-Hump Camel Function:} Contains six local minima and several saddle points.
    \begin{equation}
        f(x, y) = (4 - 2.1x^2 + \frac{x^4}{3})x^2 + xy + (-4 + 4y^2)y^2
    \end{equation}
\end{itemize}

\begin{quote}
\textit{Remark on Saddle Point Dynamics:} 
In the neighborhood of a saddle point where $||\nabla f|| \approx 0$, the update of first-order methods $\Delta \mathbf{x} = -\alpha \nabla f$ vanishes, leading to stagnation. However, the Quadratic Gradient utilizes the spectral information of the Hessian proxy. For directions where the eigenvalue $\lambda_i$ vanishes, the adaptive scaling factor $\eta_i = (\epsilon + \lambda_i)^{-1}$ compensates for the diminishing gradient. This effectively reshapes the vector field near the saddle point, ensuring a non-zero, high-velocity traversal along the directions of minimal curvature.
\end{quote}

Empirical results show that the QG variant outperforms standard GD in terms of iteration count to convergence across most benchmarks. Notably, in the Monkey Saddle experiment, the QG variant successfully avoids the stagnation observed in SGD. By identifying minimal eigenvalues $\lambda_i$, the algorithm assigns a massive update rate to those dimensions, effectively "teleporting" the trajectory away from the saddle point.

\noindent \textbf{Performance Discussion.} 

\begin{figure}[htbp]
\centering
\captionsetup[subfigure]{justification=centering}

\subfloat[The iDASH dataset]{%
    \includegraphics[width=0.48\textwidth]{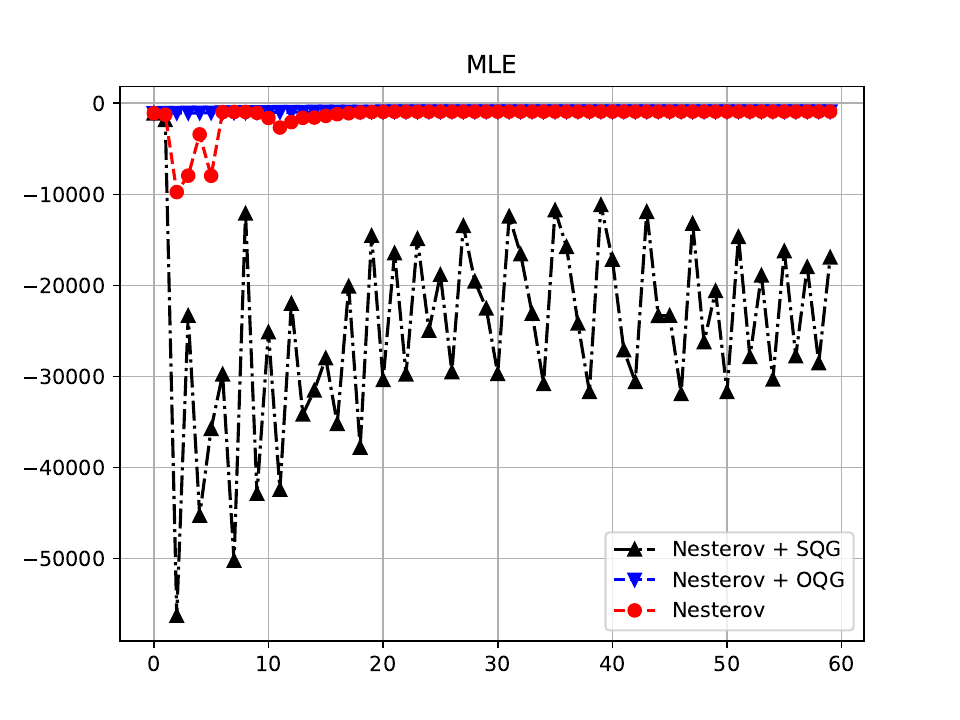}
    \label{fig:subfig01}
}
\hfill
\subfloat[The Edinburgh datasetn]{%
    \includegraphics[width=0.48\textwidth]{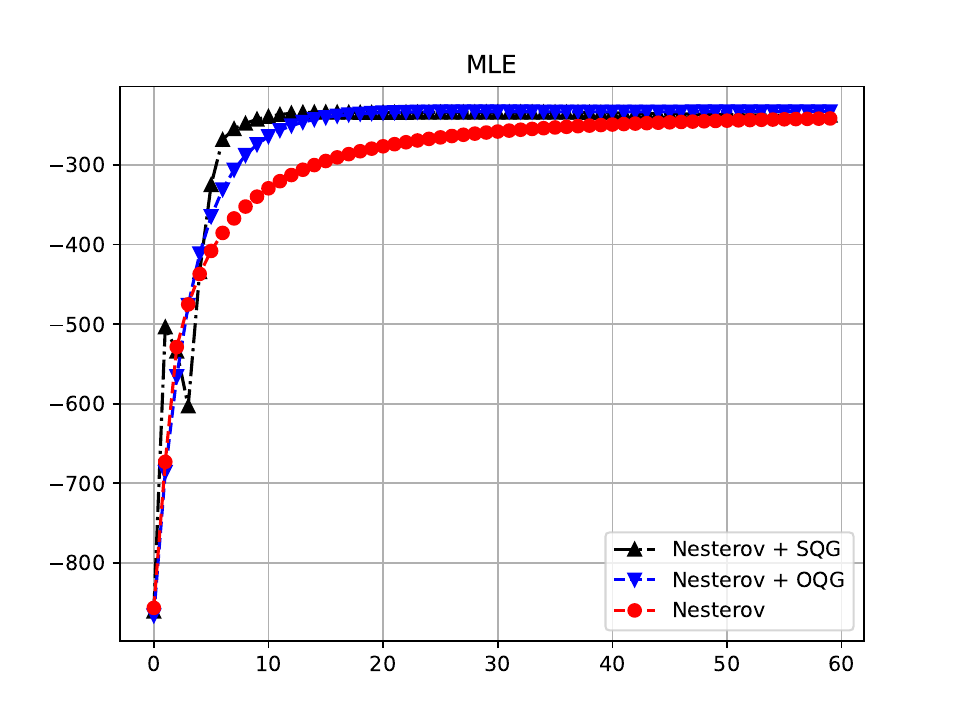}
    \label{fig:subfig02}
}

\vspace{1em} 

\subfloat[The lbw dataset]{%
    \includegraphics[width=0.48\textwidth]{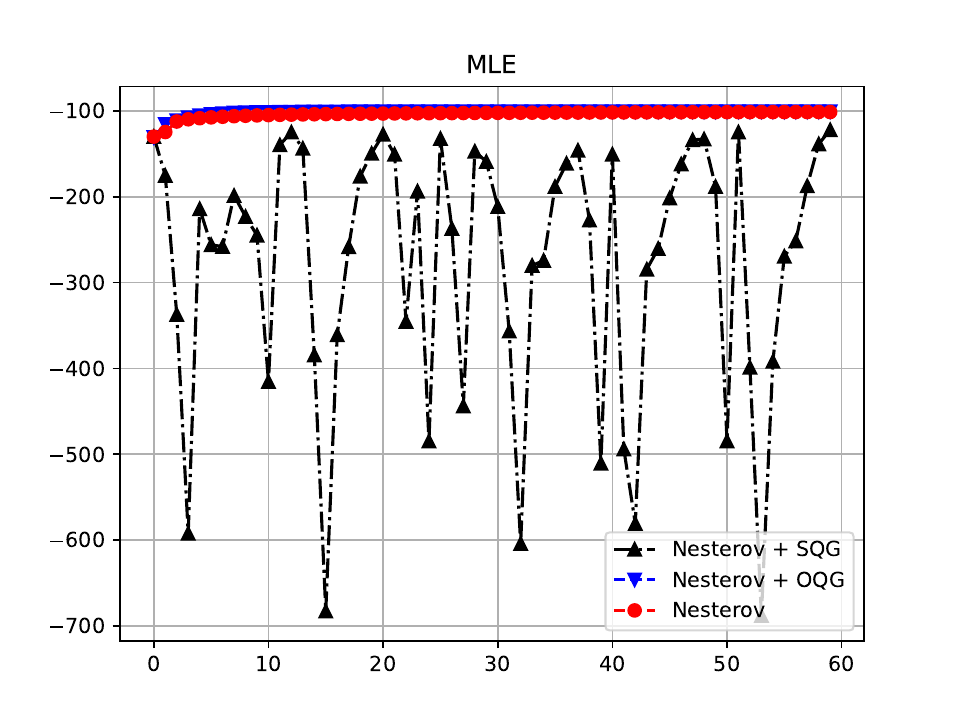}
    \label{fig:subfig03}
}
\hfill
\subfloat[The nhanes3 dataset]{%
    \includegraphics[width=0.48\textwidth]{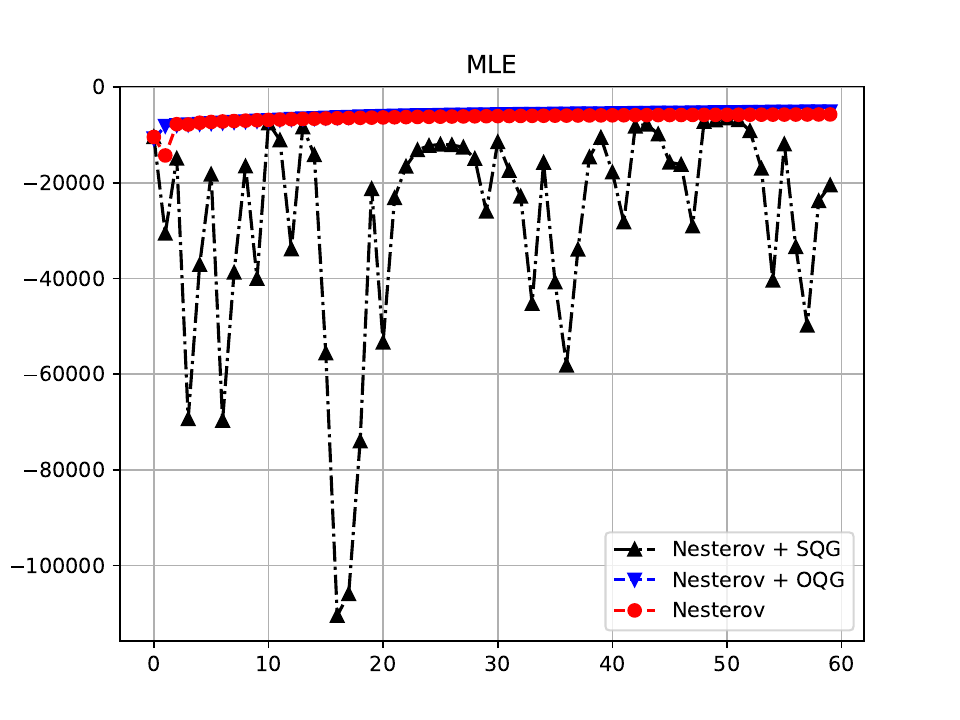}
    \label{fig:subfig04}
}

\vspace{1em} 

\subfloat[The pcs dataset]{%
    \includegraphics[width=0.48\textwidth]{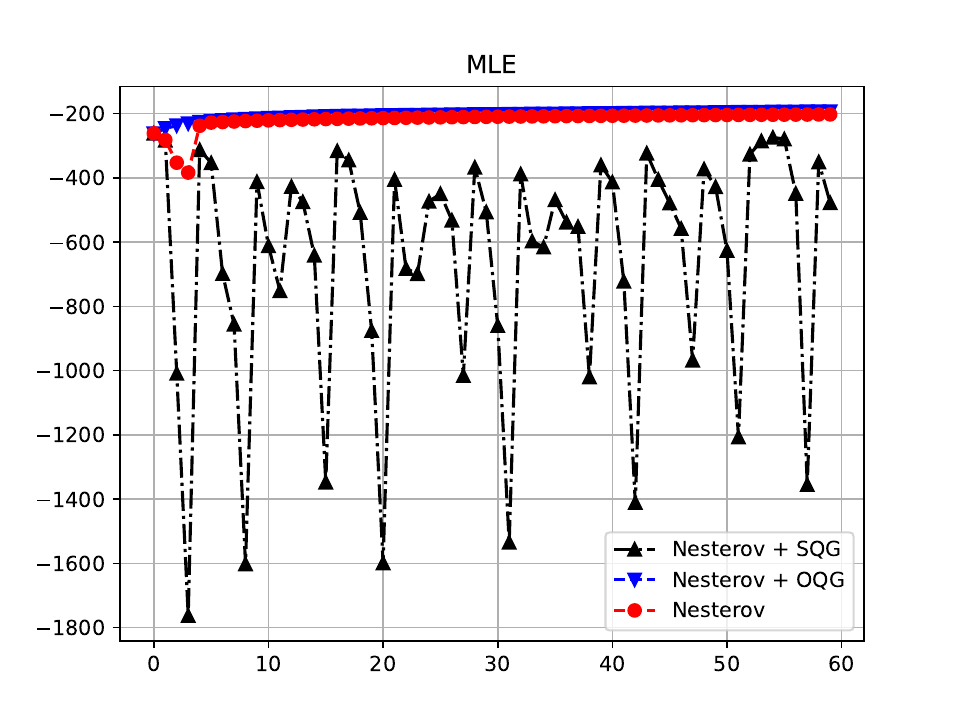}
    \label{fig:subfig03}
}
\hfill
\subfloat[The uis dataset]{%
    \includegraphics[width=0.48\textwidth]{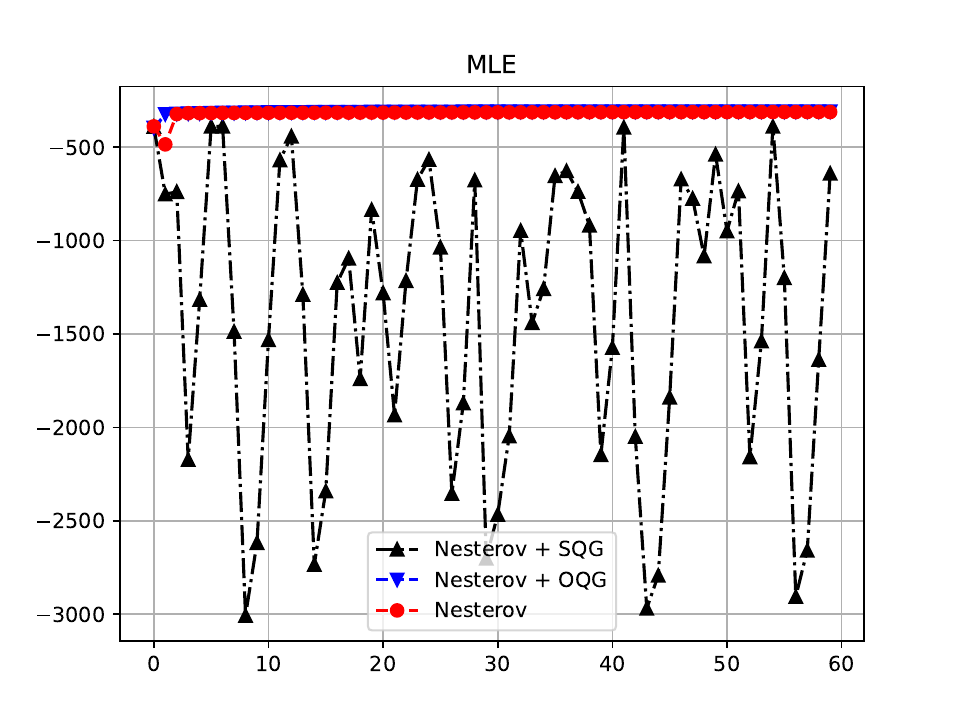}
    \label{fig:subfig04}
}

\vspace{1em} 

\subfloat[restructured MNIST dataset]{%
    \includegraphics[width=0.48\textwidth]{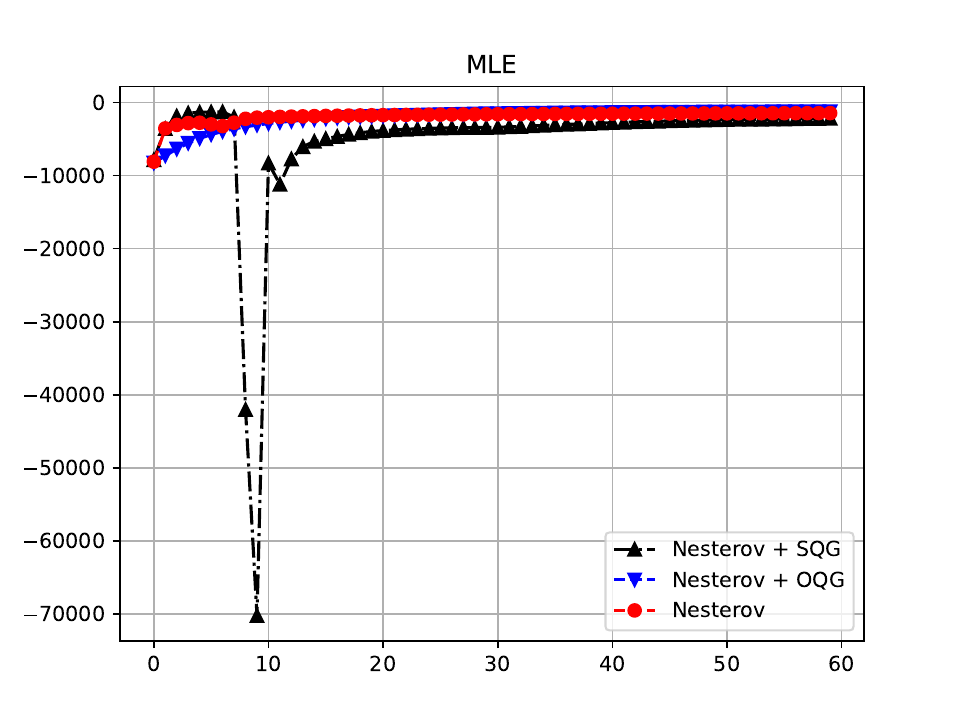}
    \label{fig:subfig03}
}
\hfill
\subfloat[The private financial dataset]{%
    \includegraphics[width=0.48\textwidth]{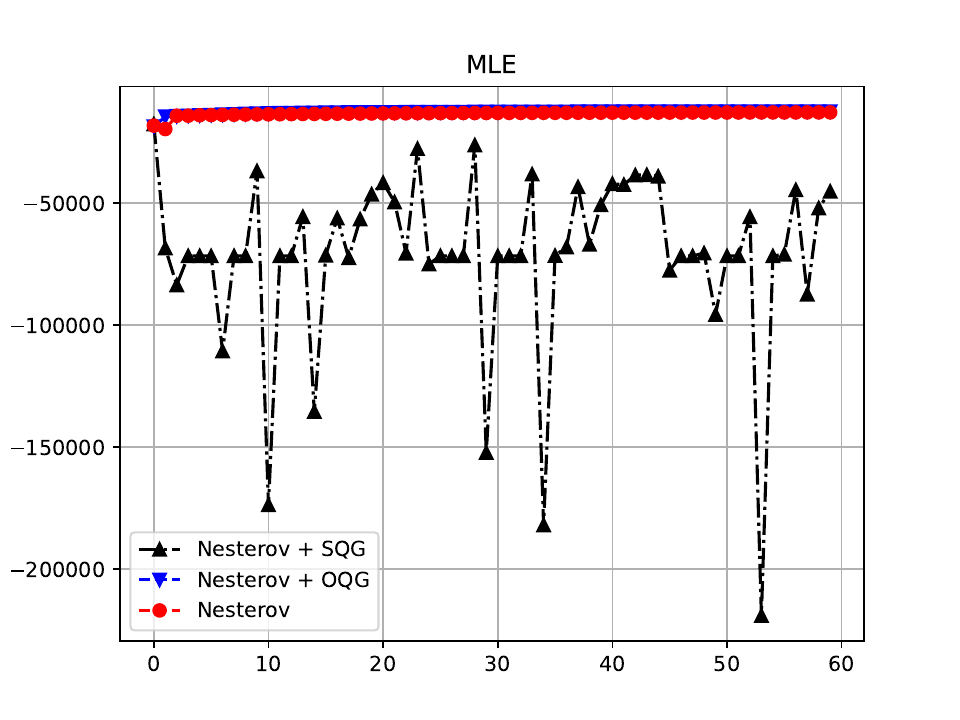}
    \label{fig:subfig04}
}

\caption{The training results of NAG + SQG vs. NAG + OQG vs. NAG in the clear domain.}
\label{fig0}
\end{figure}

\begin{figure}[htbp]
\centering
\captionsetup[subfigure]{justification=centering}

\subfloat[The iDASH dataset]{%
    \includegraphics[width=0.48\textwidth]{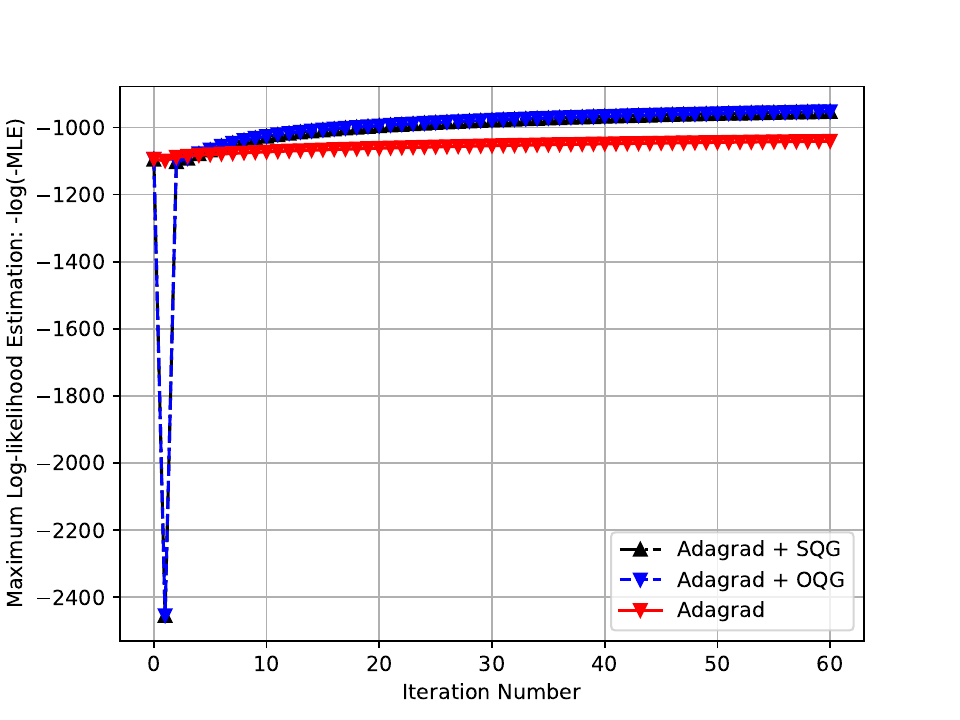}
    \label{fig:subfig01}
}
\hfill
\subfloat[The Edinburgh datasetn]{%
    \includegraphics[width=0.48\textwidth]{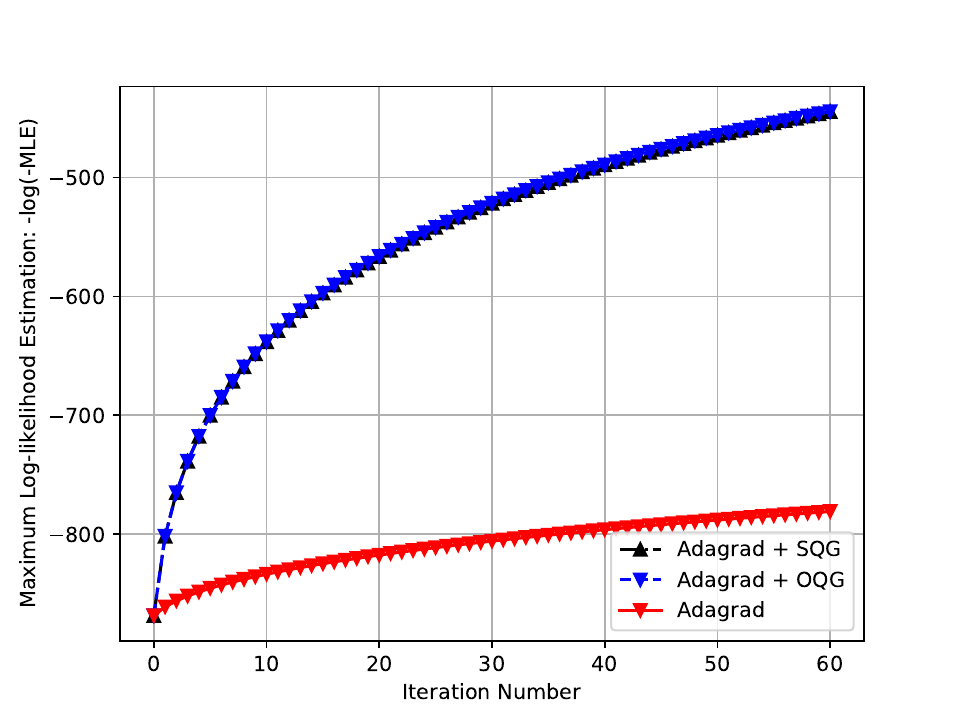}
    \label{fig:subfig02}
}

\vspace{1em} 

\subfloat[The lbw dataset]{%
    \includegraphics[width=0.48\textwidth]{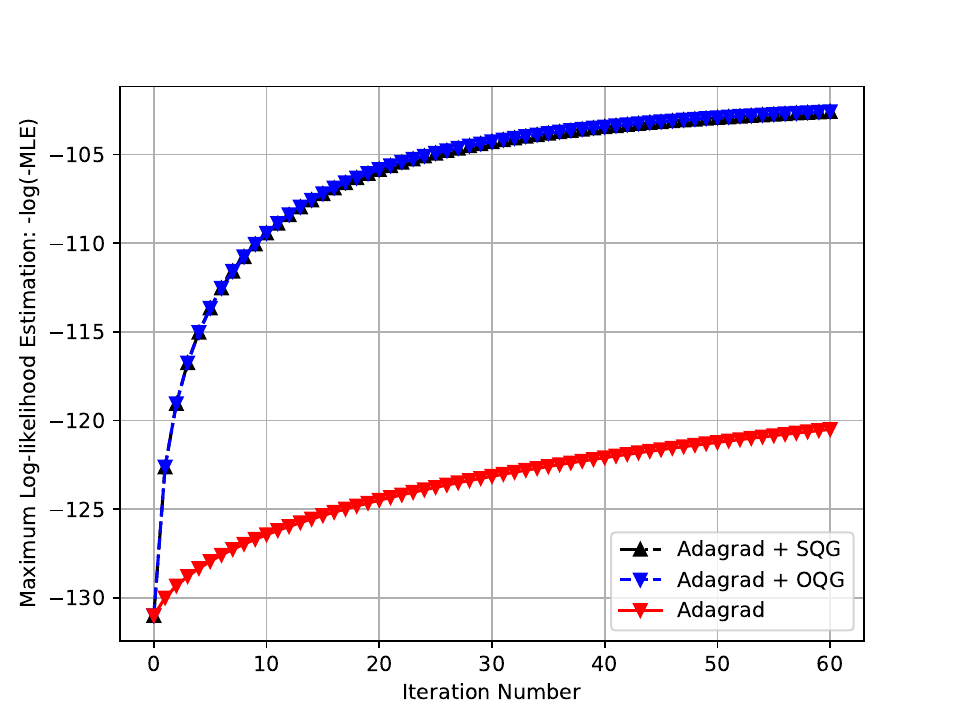}
    \label{fig:subfig03}
}
\hfill
\subfloat[The nhanes3 dataset]{%
    \includegraphics[width=0.48\textwidth]{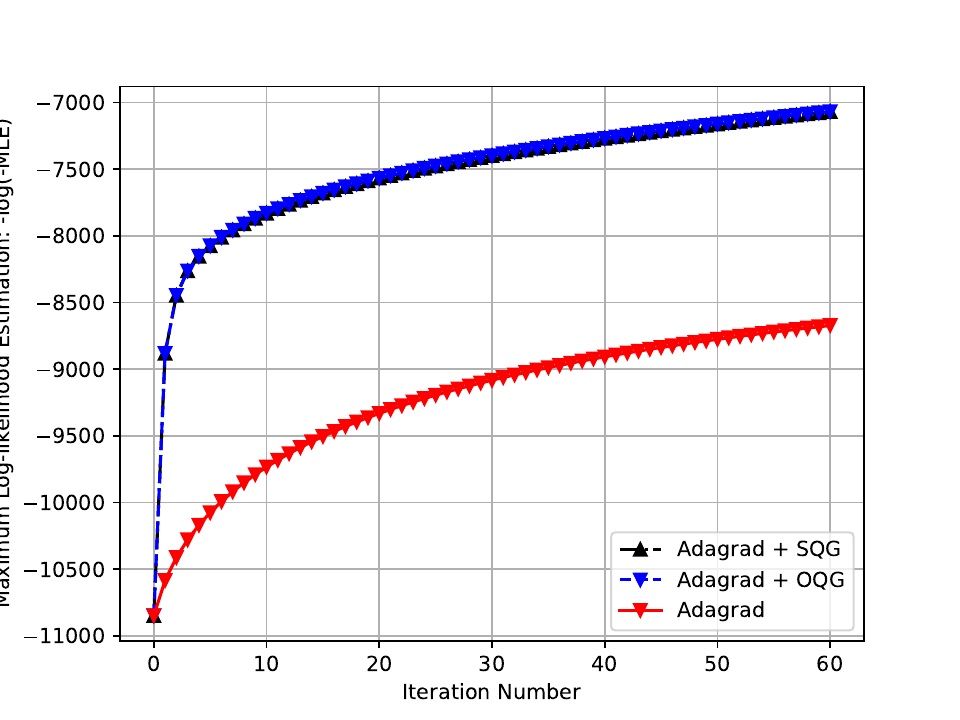}
    \label{fig:subfig04}
}

\vspace{1em} 

\subfloat[The pcs dataset]{%
    \includegraphics[width=0.48\textwidth]{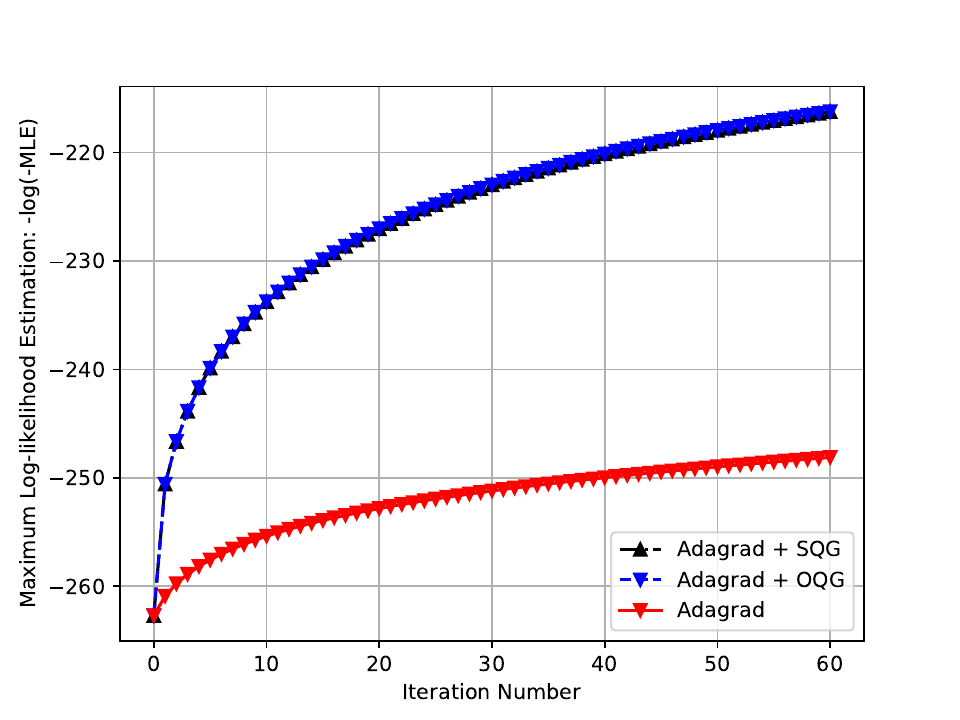}
    \label{fig:subfig03}
}
\hfill
\subfloat[The uis dataset]{%
    \includegraphics[width=0.48\textwidth]{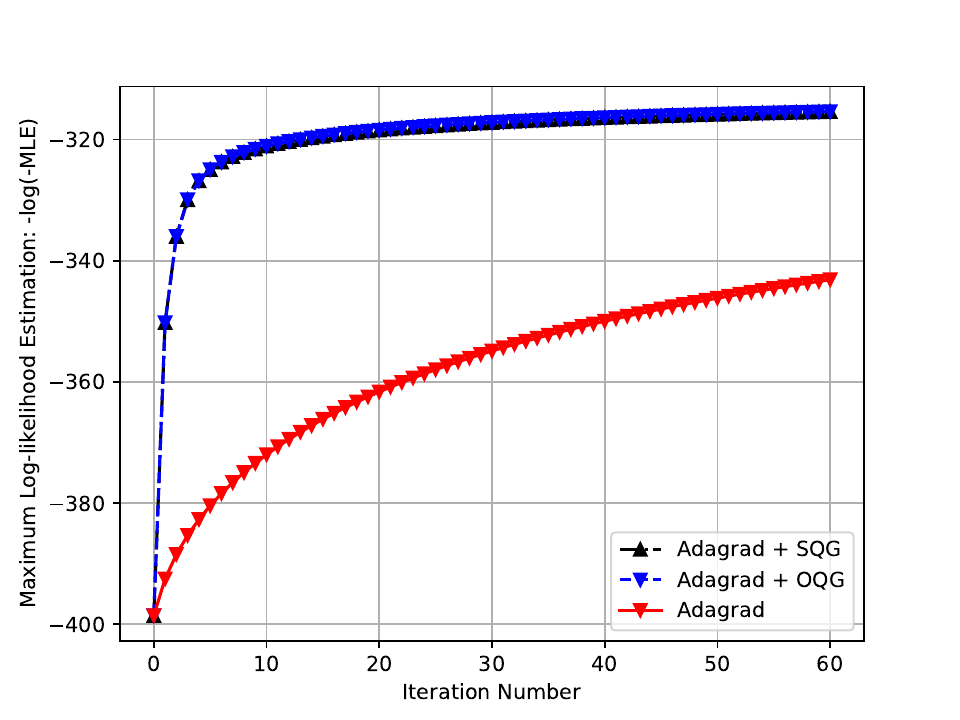}
    \label{fig:subfig04}
}

\vspace{1em} 

\subfloat[restructured MNIST dataset]{%
    \includegraphics[width=0.48\textwidth]{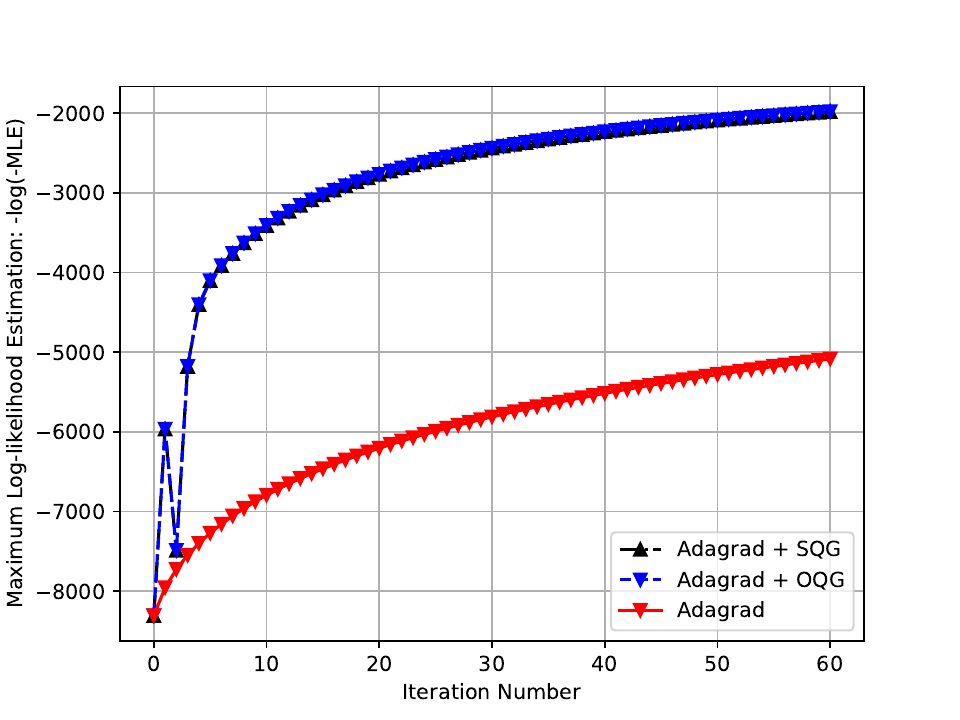}
    \label{fig:subfig03}
}
\hfill
\subfloat[The private financial dataset]{%
    \includegraphics[width=0.48\textwidth]{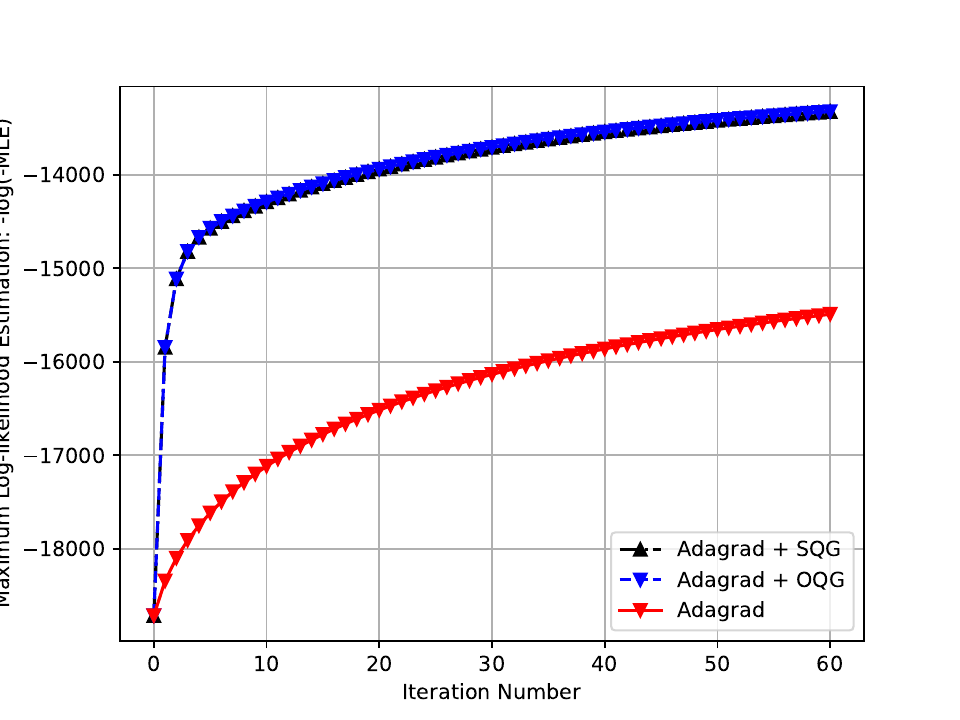}
    \label{fig:subfig04}
}

\caption{The training results of AdaGrad + SQG vs. AdaGrad + OQG vs. AdaGrad in the clear domain.}
\label{fig1}
\end{figure}

\begin{figure}[htbp]
\centering
\captionsetup[subfigure]{justification=centering}

\subfloat[The iDASH dataset]{%
    \includegraphics[width=0.48\textwidth]{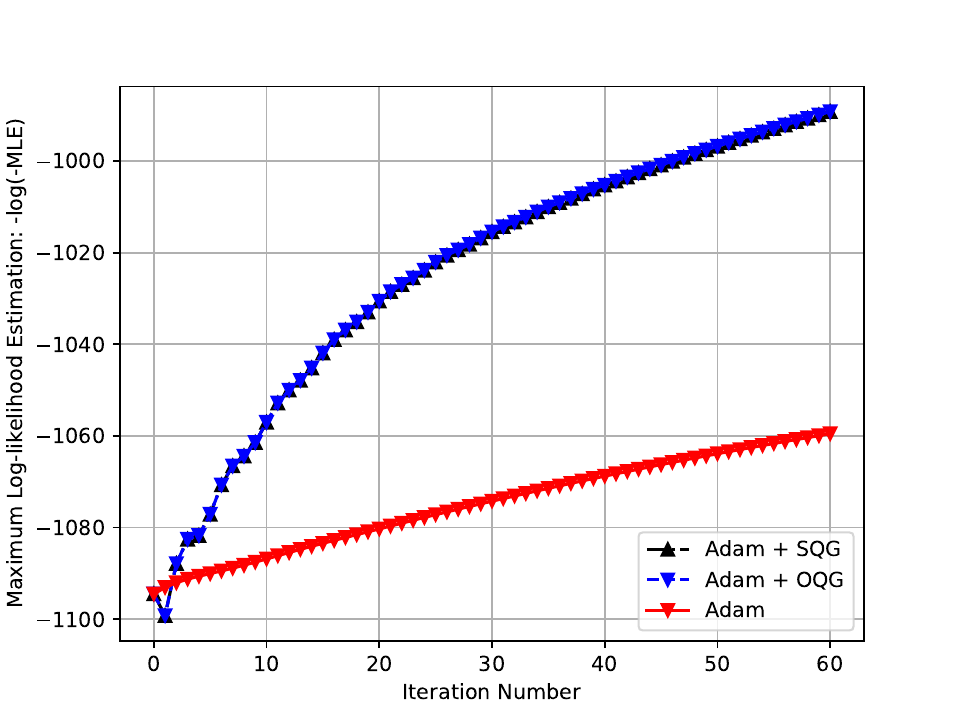}
    \label{fig:subfig01}
}
\hfill
\subfloat[The Edinburgh datasetn]{%
    \includegraphics[width=0.48\textwidth]{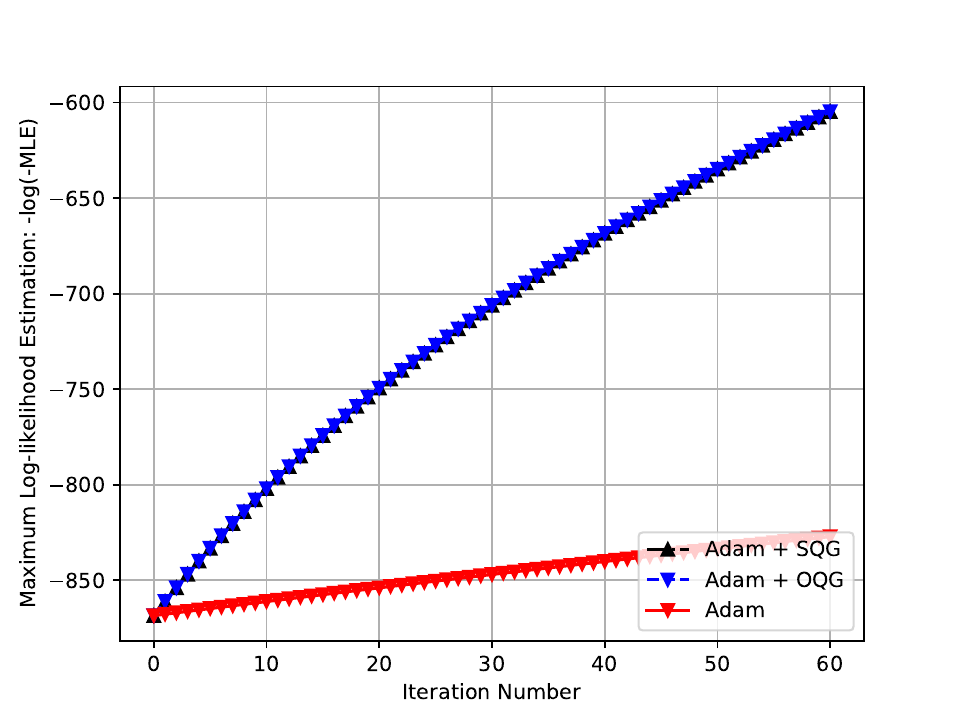}
    \label{fig:subfig02}
}

\vspace{1em} 

\subfloat[The lbw dataset]{%
    \includegraphics[width=0.48\textwidth]{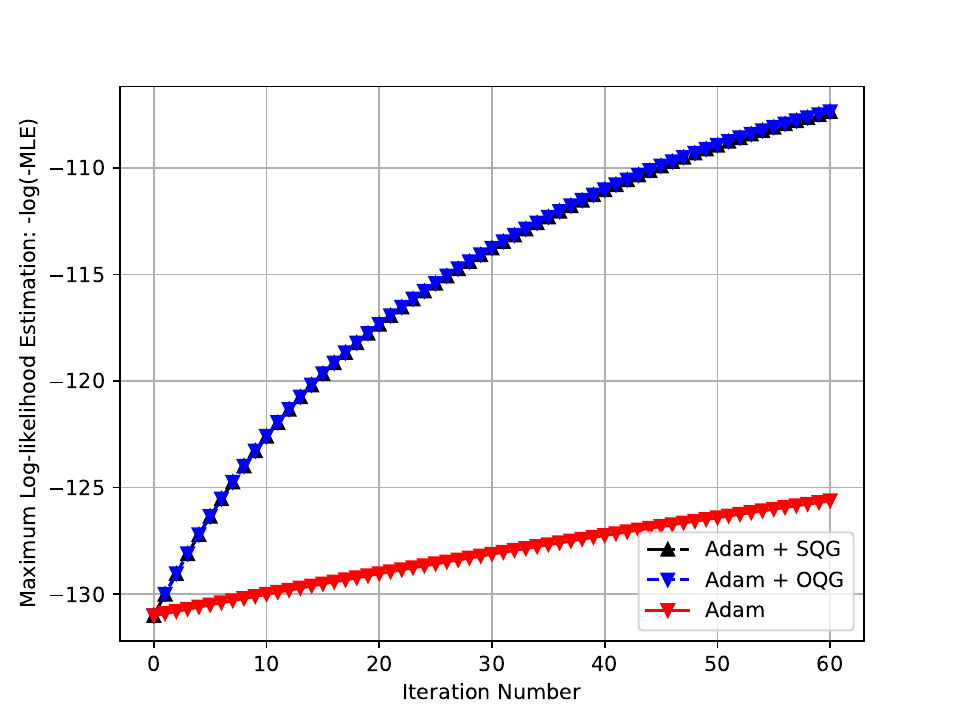}
    \label{fig:subfig03}
}
\hfill
\subfloat[The nhanes3 dataset]{%
    \includegraphics[width=0.48\textwidth]{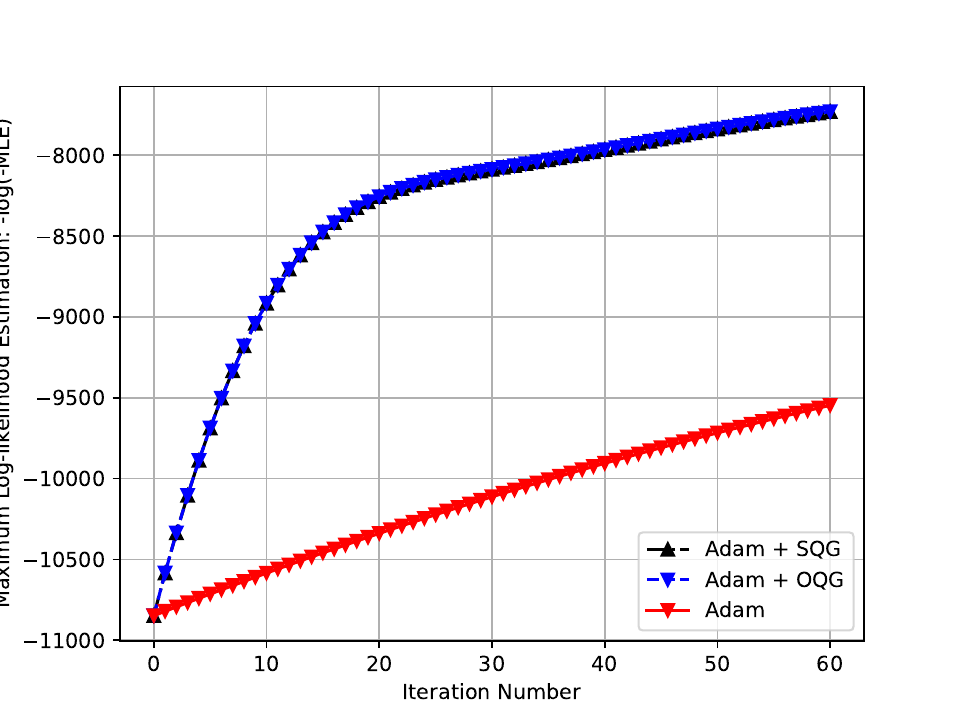}
    \label{fig:subfig04}
}

\vspace{1em} 

\subfloat[The pcs dataset]{%
    \includegraphics[width=0.48\textwidth]{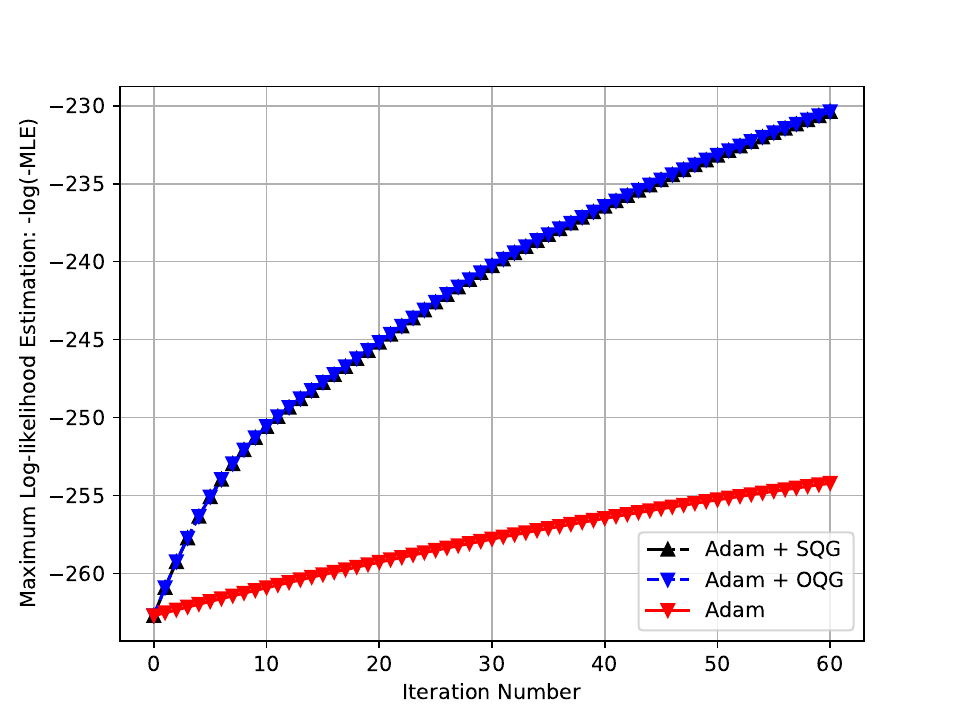}
    \label{fig:subfig03}
}
\hfill
\subfloat[The uis dataset]{%
    \includegraphics[width=0.48\textwidth]{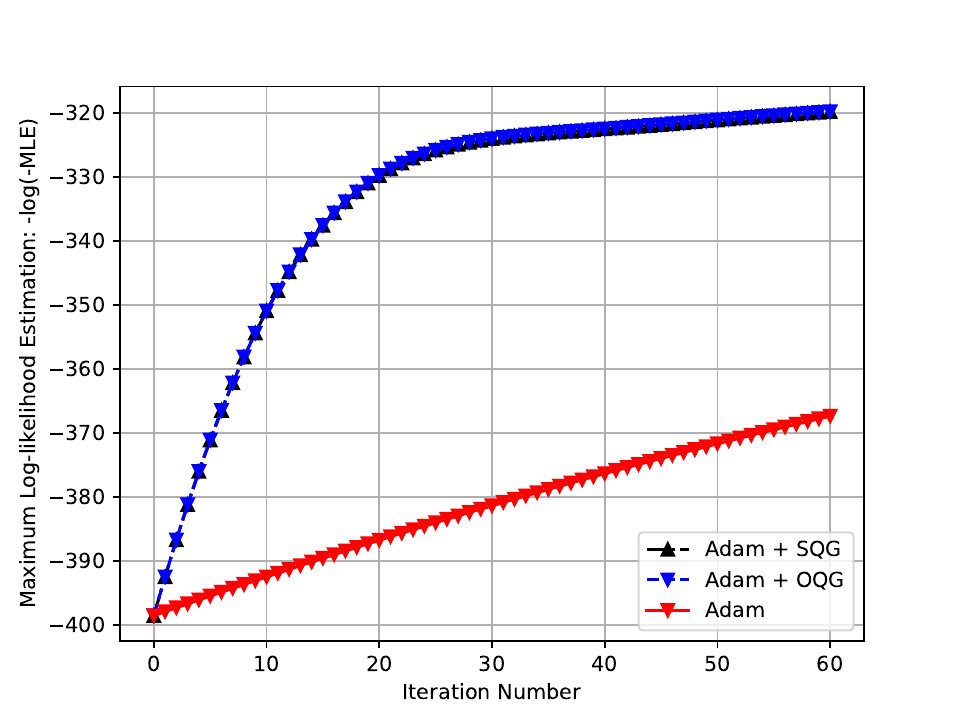}
    \label{fig:subfig04}
}

\vspace{1em} 

\subfloat[restructured MNIST dataset]{%
    \includegraphics[width=0.48\textwidth]{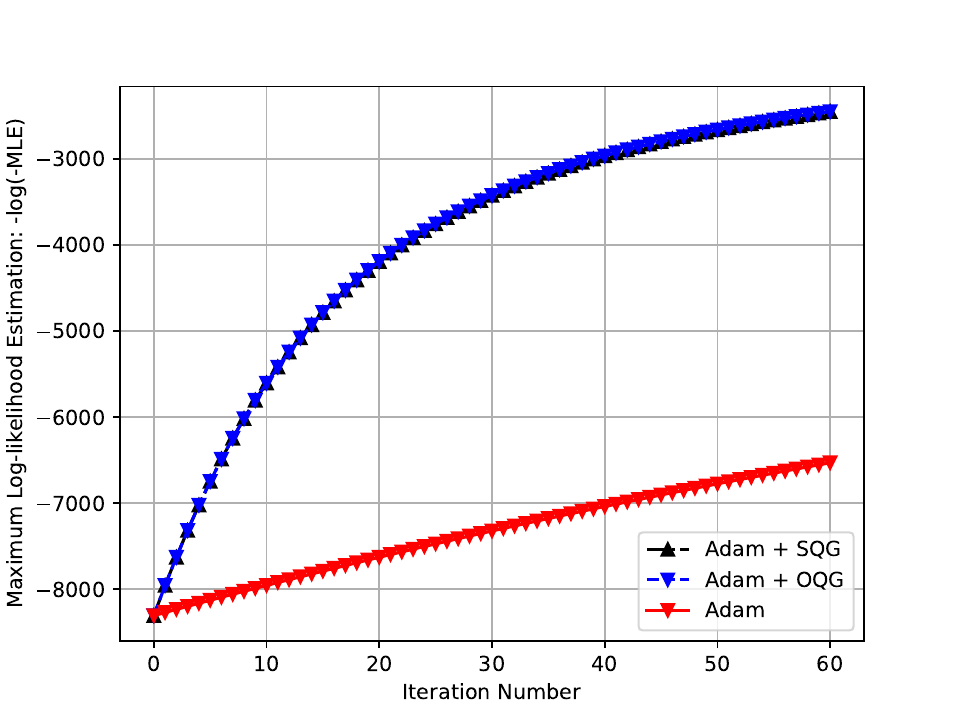}
    \label{fig:subfig03}
}
\hfill
\subfloat[The private financial dataset]{%
    \includegraphics[width=0.48\textwidth]{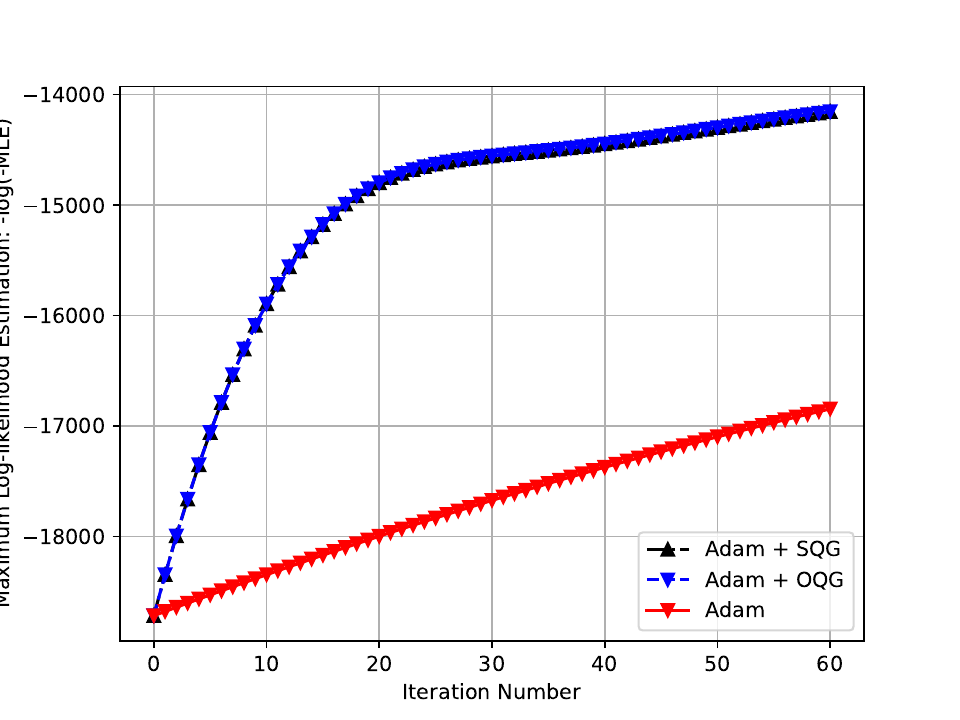}
    \label{fig:subfig04}
}

\caption{The training results of Adam + SQG vs. Adam + OQG vs. Adam in the clear domain.}
\label{fig2}
\end{figure}

As illustrated in the experimental results~\ref{fig0, fig1, fig2}, the two quadratic gradient variants for AdaGrad and Adam share a highly consistent parameter structure and exhibit remarkably similar convergence performance. Given that our Simplified Quadratic Gradient framework offers a significantly more streamlined construction process while maintaining comparable efficiency, we conclude that the SQG approach is particularly well-suited for integration with adaptive optimization methods like AdaGrad and Adam.

The observed instability and divergence in the SQG-enhanced NAG variant are not unexpected. From a theoretical perspective, the convergence of the original quadratic gradient algorithm is anchored in the Fixed Hessian framework, which permits a learning rate lower bound as high as $1$. However, our SQG framework does not strictly satisfy the stringent convergence conditions required by the Fixed Hessian method. 

Specifically, unlike the original method where the Hessian approximation provides a stable curvature estimate, the simplified nature of SQG necessitates a more conservative step-size policy for NAG-based updates. Consequently, while the original framework allows for a persistent learning rate, the SQG-enhanced NAG requires a learning rate that eventually decays below $1$ and asymptotically approaches $0$ to ensure numerical stability and global convergence.

\section{Conclusion}
In this work, we introduced the Quasi-Quadratic Gradient as a new direction for quasi-Newton optimization. By refining the search direction through the interaction of the inverse Hessian proxy and the gradient, we achieve faster convergence in complex optimization landscapes. Future work will explore the application of QQG in large-scale distributed training.

All the python source code to implement the experiments in the paper  is openly available at: \href{https://github.com/petitioner/ML.QuasiQuadraticGradient}{$\texttt{https://github.com/petitioner/ML.QuasiQuadraticGradient}$}  .

\bibliography{ML.QuasiQuadraticGradient}
\bibliographystyle{apalike}

\end{document}